\documentclass[12pt, final]{article}
\usepackage{a4}
\usepackage{amsmath}%
\usepackage{amstext}%
\usepackage{amssymb}%
\usepackage{showkeys}%
\usepackage{epsfig}%
\usepackage{cite}
\newtheorem{theorem}{Theorem}

\newtheorem{proposition}[theorem]{Proposition}
\newtheorem{remark}{Remark}

\newcounter{unnumber}

\newcommand{\R}{\mathbb{R}}%
\newcommand{\N}{\mathbb{N}}%
\DeclareMathOperator*\cl{cl}%
\DeclareMathOperator*\dom{dom}%

\DeclareMathOperator*\core{core}%
\DeclareMathOperator*\co{co}%
\DeclareMathOperator*\B{\overline{\R}}%
\DeclareMathOperator*\im{Im}%
\DeclareMathOperator*\pr{pr}%
\DeclareMathOperator*\ba{bar}%

\title{On extension results for $n$-cyclically monotone operators in reflexive Banach spaces}

 \author{Radu Ioan Bo\c{t} \thanks
 {Faculty of Mathematics, Chemnitz University of Technology,
D-09107 Chemnitz, Germany, e-mail:
 radu.bot@mathematik.tu-chemnitz.de. Research partially supported by DFG (German Research Foundation), project WA 922/1-3.} \and Ern\"{o} Robert Csetnek
 \thanks {Faculty of Mathematics, Chemnitz University of Technology,
D-09107 Chemnitz, Germany, e-mail:
 robert.csetnek@mathematik.tu-chemnitz.de}}

\begin{document}
\maketitle

\noindent \textbf{Abstract.} In this paper we provide some extension results for $n$-cyclically monotone operators in reflexive Banach spaces
by making use of the Fenchel duality. In this way we give a positive answer to a question posed by Bauschke and Wang in \cite{bauschke-wang}.\\



\section{Introduction and preliminaries}

This paper is motivated by the work of Bauschke and Wang \cite{bauschke-wang}, where extension results for $n$-cyclically monotone operators in Hilbert spaces are delivered. In that paper the authors provide a new proof of the main result in Voisei's article \cite{voiseiext}, where refinements of the Debrunner-Flor theorem (cf. \cite{debrunner-flor}) for $n$-cyclically monotone operators are obtained, for the proof of which techniques relying on fixed point results are used. Different to the approaches in \cite{voiseiext}, Bauschke and Wang make use of the \textit{Fitzpatrick function} associated with a $n$-cyclically monotone operator, well-studied in \cite{bbbrw}, and of the convex duality theory. The Fitzpatrick function for a $n$-cyclically monotone operator has been introduced
and investigated in \cite{bbbrw} as an extension of the one considered by Fitzpatrick in \cite{fitz} for the study of monotone operators and which played in the last years an important role in the development of this field.

Since the main result in \cite{bauschke-wang} is stated in Hilbert spaces and its proof uses in a determinant manner the characteristics of this framework,
the authors of the paper ask in \cite[Remark 3.5 (6)]{bauschke-wang} whether or not is it possible to extend the result to Banach spaces. In the following we give a positive answer to this question in the setting of reflexive Banach spaces. Actually we are able to formulate and prove three extension results, differing in the hypotheses assumed.

The structure of this work is the following. In the next two subsections of the first section we introduce some elements of convex analysis as well as the notion of a $n$-cyclically monotone operator along with some of its properties, respectively. The second section is dedicated to the formulation of the extension results, while in the third one we formulate some conclusions and propose some possible further research.

\subsection{Elements of convex analysis}

We start by considering a real separated locally convex space $X$ and its continuous dual space $X^*$. The notation $\langle x^*,x\rangle$ stands for the value of the linear continuous functional $x^*\in X^*$ at $x\in
X$. The notation $p=\langle\cdot,\cdot\rangle$ is used for the pairing between $X^*$ and $X$. For a subset $C$ of $X$ we denote by $\cl C, \co C$ and $\core C$ its \emph{closure}, \emph{convex hull} and \emph{algebraic
interior} (or \emph{core}), respectively. Note that if $C$
is a convex set, then an element $x\in X$ belongs to $\core C$ if
and only if $\bigcup_{\lambda\geq 0}\lambda(C-x)=X$ (see also \cite{Rock-conj-dual, Zal-carte}).

For a function $f:X\rightarrow\B=\R\cup\{\pm\infty\}$ we denote by
$\dom f=\{x\in X:f(x)<+\infty\}$ its \emph{domain}.
We call $f$ \emph{proper} if $\dom f\neq\emptyset$ and
$f(x)>-\infty$ for all $x\in X$. For a function $f:A\times B\rightarrow\B$, where $A$ and $B$ are
nonempty sets, we denote by $f^\top$ the \emph{transpose} of $f$,
namely the function $f^\top:B\times A\rightarrow\B,
f^\top(b,a)=f(a,b)$ for all $(b,a)\in B\times A$. Here we also introduce the \emph{projection
operator} $\pr_A:A\times B\rightarrow A$, $\pr_A(a,b)=a$ for all
$(a,b)\in A\times B$.

The \emph{Fenchel-Moreau (Legendre-Fenchel) conjugate} of $f$ is the function
$f^*:X^*\rightarrow\B$ defined by $$f^*(x^*)=\sup\limits_{x\in
X}\{\langle x^*,x\rangle-f(x)\} \ \forall x^*\in X^*.$$ We mention here
some important properties of a conjugate function. First of all, we have the
so-called \emph{Young-Fenchel inequality}
$$f^*(x^*)+f(x)\geq\langle x^*,x\rangle\mbox{ for all }(x,x^*)\in
X\times X^*.$$ If $f$ is proper, then $f$ is convex and lower
semicontinuous if and only if $f^{**}=f$ (see
\cite{EkTem,Zal-carte}).

For $C \subseteq X$ a given set we denote by $\delta_C : X \rightarrow \overline \R$, defined by
$$\delta_C(x) = \left \{\begin{array}{ll}
0, & \ \mbox{if} \ x \in C,\\
+\infty, & \ \mbox{otherwise}
\end{array}\right.$$
its \emph{indicator function} and by $\sigma_C : X^* \rightarrow \overline \R$, defined by $\sigma_C(x^*) = \delta^*_C(x^*)$, its \emph{support function}. The \emph{barrier cone}
of $C$ is the set $\ba C := \dom \sigma_C$ and one has $\ba C = \ba \cl \co C$.

Given a linear continuous operator $A:X\rightarrow Y$ ($Y$ is another separated locally convex space), we denote by $\im(A)$ its
\emph{image-set} $\im(A)=\{Ax:x\in X\}$ and by $A^*$ its
\emph{adjoint operator} $A^*:Y^*\rightarrow X^*$, given by
$\langle A^*y^*,x\rangle=\langle y^*,Ax\rangle \ \forall y^*\in
Y^* \ \forall x\in X$. For a function $h:X\rightarrow Y$ and a set $D \subseteq Y$ we use the notation
$h^{-1}(D)=\{x\in X:h(x)\in D\}$. Having $f,g:X\rightarrow\B$ two functions we consider
also their \emph{infimal convolution}, which is the function denoted by $f\Box
g:X\rightarrow\B$, $f\Box g(x)=\inf_{u\in X}\{f(u)+g(x-u)\}$ for
all $x\in X$. We say that the infimal convolution is \emph{exact
at $x\in X$} if the infimum in its definition is attained.
Moreover, $f\Box g$ is said to be \emph{exact} if it is exact at
every $x\in X$. When an infimum or a supremum is attained we write min and max instead
of inf, respectively, sup.

Let us recall in the following the classical Fenchel duality result to which we will refer several times in the next section.

\begin{theorem}\label{dual-f} (Fenchel-Rockafellar duality, \cite{Rock-conj-dual}, \cite[Corollary 2.8.5]{Zal-carte}) Let $X$ and $Y$ be separated locally convex spaces, $A:X\rightarrow Y$ a linear and continuous operator and $f:X\rightarrow \B$ and $g:Y\rightarrow\B$ be two proper and convex functions such that one of the following regularity conditions is fulfilled: \begin{center}
\begin{tabular}{r|l}
$(RC_1)$ \ & \ $\exists x_0\in\dom f\cap A^{-1}(\dom g)$ such that $g$ is continuous at $Ax_0$;
\end{tabular}
\end{center}

\begin{center}
\begin{tabular}{r|l}
$(RC_2)$ \ & \ $X,Y$ are Fr\'{e}chet spaces, $f$ and $g$ are lower semicontinuous and \\ \ & \ $0\in\core\big(\dom g-A(\dom f)\big)$.\end{tabular}
\end{center}

\noindent Then $$\inf_{x\in X}\{f(x)+g(Ax)\}=\max_{y^*\in Y^*}\{-f^*(-A^*y^*)-g^*(y^*)\}.$$
\end{theorem}

\begin{remark}\rm Let us notice that instead of the core, one can use in the above duality result other generalized interiority notions, like the \emph{intrinsic core}, or the \emph{strong quasi relative interior}. We refer to \cite{bgw-carte, Zal-art, Zal-carte, Gowda-Teboulle, b-hab} for further considerations concerning generalized interior-type regularity conditions ensuring the above duality result. We remark that in case $X,Y$ are Fr\'{e}chet spaces, $f$ and $g$ are proper, convex and lower semicontinuous then $(RC_1)\Rightarrow (RC_2)$.
\end{remark}

Consider in the following that $(X,\|\cdot\|)$ is a real normed space. We say that $f:X\rightarrow \B$ is \emph{coercive} if $\lim_{\|x\|\rightarrow+\infty}f(x)=+\infty$. It is obvious that $f$ is coercive if and only if all level sets $[f\leq\lambda]:=\{x\in X:f(x)\leq\lambda\}$, $\lambda\in\R$, are bounded. It follows by \cite[Theorem 7A(a)]{rock-level-cont} that if $f$ is a proper, convex, lower semicontinuous and coercive function, then $f^*$ is finite and continuous at $0$ (see also \cite[Exercise 2.41]{Zal-carte}). The function $f:X\rightarrow \B$ is said to be \emph{strongly coercive} if $\lim_{\|x\|\rightarrow+\infty}f(x)/\|x\|=+\infty$. In view of \cite[Lemma 3.6.1]{Zal-carte}, if $f$ is a proper, convex, lower semicontinuous and strongly coercive function, then $\dom f^*=X^*$. In this case $f^*$ is continuous on $X^*$ with respect to the strong topology. This is a direct consequence of \cite[Corollary 2.5]{EkTem}, by noticing that the function $f^*$ is lower semicontinuous with respect to the strong topology of $X^*$, since it is lower semicontinuous with respect to the weak$^*$ topology of $X^*$.

Having a linear and continuous operator $B:X \rightarrow X^*$ we call it \emph{coercive} (\emph{strongly coercive}) if $x \mapsto \langle Bx,x\rangle$ is a coercive (strongly coercive) function.

\subsection{$n$-cyclically monotone operators}

We recall in this subsection some basic facts regarding $n$-cyclically monotone operators. Consider in the following a real Banach space $(X,\|\cdot\|)$ with corresponding dual space $X^*$. For a set-valued operator $S:X\rightrightarrows X^*$ we use the notations $G(S):=\{(x,x^*)\in X\times X^*:x^*\in S(x)\}$, $D(S):={\pr}_X G(S)=\{x\in X:S(x)\neq\emptyset\}$ and $R(S):={\pr}_{X^*}G(S)=\cup\{S(x):x\in D(S)\}$ for its \emph{graph}, \emph{domain}, respectively, \emph{range}. The operator $S$ is said to be \emph{$n$-monotone} (or \emph{$n$-cyclically monotone}), where $n\in\N,n\geq2$, if $$\sum\limits_{i=1}^{n}\langle s_i^*,s_{i+1}-s_i\rangle\leq 0 \ \forall (s_i,s_i^*)\in G(S) \ \mbox{ with } s_{n+1}=s_1.$$ Let us notice that $2$-monotonicity is nothing else than the classical \textit{monotonicity}, that is $\langle x^*-y^*,x-y\rangle\geq 0$ for all $(x,x^*),(y,y^*)\in G(S)$. The operator $S$ is \emph{cyclically monotone} if $S$ is $n$-cyclically monotone for all $n\in\{2,3,...\}$. The multifunction $S$ is \emph{maximal $n$-monotone} if $S$ is $n$-monotone and no proper extension (in the sense of inclusion of graphs) of $S$ is $n$-monotone. Let us note that $S$ is maximal $2$-monotone exactly when $S$ is \emph{maximal monotone} (we refer to \cite{phelps, simons} for more on this classical notion). One of the important results concerning cyclically monotone operators is due to Rockafellar, who proved in the finite dimensional setting that maximal cyclically monotone operators are exactly the subdifferential operators of proper, convex and lower semicontinuous functions (cf. \cite[Theorem 24.9]{Rock-carte}).

Further, let us consider the \emph{Fitzpatrick function of order n} associated with $S$ (cf. \cite{bbbrw}), $F_{S,n}:X\times X^*\rightarrow \B$, which plays an significant role in the next section: $$F_{S,n}(x,x^*)=\sup\limits_{\substack{(s_i,s_i^*)\in G(S)\\i=\overline{1,n-1}}}\left\{\sum_{i=1}^{n-2}\langle s_i^*,s_{i+1}-s_i\rangle+ \langle s_{n-1}^*, x-s_{n-1}\rangle+\langle x^*,s_1\rangle\right\}.$$ For $n=2$ we obtain the classical \emph{Fitzpatrick function} introduced and investigated in \cite{fitz}, $F_{S,2}:X\times X^*\rightarrow\B$, $$F_{S,2}(x,x^*)=\sup\{\langle x^*,s\rangle+\langle s^*,x\rangle-\langle s^*,s\rangle:(s,s^*)\in G(S)\}.$$ The Fitzpatrick function plays an indisputable role in the modern monotone operator theory due to the fact that it links the duality results in convex analysis with the property of maximal monotonicity for operators. We refer to \cite{bbbrw, bausch-m-s, bauschke-wang, BCW-op-nonlinear, BGW-max-comp, b, bu-sv-02, mas, penot-zal, simons, voisei} for more details concerning this fact.

Let us recall in the following some results regarding $n$-monotone operators.

\begin{proposition}\label{ext-reun} (cf. \cite[Proposition 2.7]{bbbrw}) Let $S:X\rightrightarrows X^*$ be $n$-monotone for some $n\in\{2,3,...\}$, $(x,x^*)\in X\times X^*$ and let us define $T:X\rightrightarrows X^*$ via $G(T)=G(S)\cup\{(x,x^*)\}$. Then \begin{equation}\label{ext-reun-rel}T\mbox{ is $n$-monotone}\Leftrightarrow F_{S,n}(x,x^*)\leq\langle x^*,x\rangle.
\end{equation}\end{proposition}

The following result was proved in \cite{bauschke-wang} in the setting of Hilbert spaces. By using the same techniques one can show that it remains valid in the framework of reflexive Banach spaces.

\begin{proposition}\label{prop-dom} (cf. \cite[Proposition 2.6]{bauschke-wang}) Let $X$ be a reflexive Banach space and $S:X\rightrightarrows X^*$ be a given multifunction. Then \begin{equation}\label{dom-2}\co G(S)\subseteq\dom F_{S,2}^{*\top}\subseteq \cl\co G(S)\subseteq \cl\co D(S)\times\cl\co R(S)\end{equation} and \begin{equation}\label{dom-n}\forall n\in\{3,4,...\} \ \co D(S)\times\co R(S)\subseteq\dom F_{S,n}^{*\top}\subseteq \cl\co D(S)\times\cl\co R(S).\end{equation}
\end{proposition}

\begin{remark}\label{s-tr}\rm (see also \cite[Remark 2.9]{bauschke-wang}) Let $S:X\rightrightarrows X^*$ be a multifunction, $w^* \in X^*$ and define $S':X\rightrightarrows X^*$ by $S'(x)=-w^*+S(x)$ for all $x\in X$. One can prove that $S'$ is $n$-monotone if and only if $S$ is $n$-monotone. Further, $F_{S',n}(x,x^*)=F_{S,n}(x,x^*+w^*)-\langle w^*,x\rangle$ and $F^*_{S',n}(x^*,x)=F^*_{S,n}(x^*+w^*,x)-\langle w^*,x\rangle$ for all $(x,x^*)\in X\times X^*$. Hence, $p\leq F^*_{S',n}\Leftrightarrow p\leq F^*_{S,n}$.
\end{remark}

\section{Extension results}

We extend in this section to the setting of reflexive Banach spaces the convex-analytical approach used by Bauschke and Wang in \cite{bauschke-wang} for obtaining extension results for $n$-monotone operators.

Throughout this section $X$ is a reflexive Banach space.

\begin{theorem}\label{ext-1-coerc} Let $S:X\rightrightarrows X^*$ be $n$-monotone for some $n\in\{2,3,...\}$. Suppose that $G(S)\neq\emptyset$ and \begin{equation}\label{pleq}p\leq F^*_{S,n}.\end{equation} Consider a linear, continuous, monotone and strongly coercive operator $B:X\rightarrow X^*$. Then for every $w^*\in X^*$ there exists $x\in\cl\co D(S)$ such that $\{(x,w^*-Bx)\}\cup G(S)$ is $n$-monotone.\end{theorem}

\noindent {\bf Proof.} We denote $C:=\cl\co D(S)$. We prove first the result in case $w^*=0$. Applying Proposition \ref{ext-reun}, it is enough to show that there exists $x\in C$ such that $F_{S,n}(x,-Bx)\leq\langle -Bx,x \rangle$, or, equivalently $$\min_{x\in X}\{F_{S,n}(x,-Bx)+\langle Bx,x \rangle+\delta_C(x)\}\leq 0,$$ which is nothing else than \begin{equation}\label{dual-pb} \max_{x\in X}\{-F_{S,n}(Ax)-f(x)\}\geq 0,\end{equation} where $A:X\rightarrow X\times X^*$ and $f:X\rightarrow\B$ are defined by $Ax=(x,-Bx)$ and $f(x)=\langle Bx,x\rangle+\delta_C(x)$ for all $x\in X$, respectively. Obviously, $A$ is a linear and continuous operator and its adjoint operator is $A^*:X^*\times X\rightarrow X^*$, $A^*(x^*,x)=x^*-B^*x$ for all $(x^*,x)\in X^*\times X$. In the hypotheses we work, the function $f$ is proper, convex, lower semicontinuous and strongly coercive. This means that $f^*$ is continuous on $X^*$ with respect to its strong topology (see subsection 1.1). By using the duality result Theorem \ref{dual-f} (notice that $(RC_1)$ is fulfilled) and taking into account that $X$ is reflexive, we get that \begin{equation}\label{dualitate} \inf_{(x^*,x)\in X^*\times X}\{F_{S,n}^*(x^*,x)+f^*(-A^*(x^*,x))\}=\max_{x\in X}\{-F_{S,n}(Ax)-f(x)\},\end{equation} hence, in order to show that \eqref{dual-pb} holds, we only have to prove that \begin{equation}\label{dearatat} \inf_{(x,x^*)\in X\times X^*}\{F_{S,n}^*(x^*,x)+f^*(B^*x-x^*)\}\geq 0.\end{equation} In view of Proposition \ref{prop-dom} it remains to show that \begin{equation}\label{dearatat2} F_{S,n}^*(x^*,x)+f^*(B^*x-x^*)\geq 0 \ \forall (x,x^*)\in C\times X^*.\end{equation}

Take an arbitrary $(x,x^*)\in C\times X^*$. Then it holds (cf. \eqref{pleq})
\begin{eqnarray*} F_{S,n}^*(x^*,x)+f^*(B^*x-x^*) & \geq & \langle x^*,x\rangle + \sup_{y \in C} \{\langle B^*x-x^*,y\rangle-\langle By,y \rangle\}\\
& \geq & \langle x^*,x\rangle + \langle B^*x-x^*,x\rangle - \langle Bx,x \rangle \\ & = &  0.\end{eqnarray*}
Since $(x,x^*)\in C\times X^*$ is arbitrary chosen, the inequality \eqref{dearatat2} is fulfilled and the conclusion holds for $w^*=0$.

Assume now that $w^*\in X^*$ is arbitrary. Consider the operator $S':X\rightrightarrows X^*$ defined by $S'(x)=-w^*+S(x)$ for all $x\in X$. By Remark \ref{s-tr}, the inequality \eqref{pleq} holds for $S'$, too. Since $D(S')=D(S)$, the above considerations provide a point $(x,-Bx)\in C\times X^*$ such that $\{(x,-Bx)\}\cup G(S')$ is $n$-monotone, which is nothing else than $\{(x,w^*-Bx)\}\cup G(S)$ is $n$-monotone.\hfill{$\Box$}

\begin{remark}\label{caz-part-bauschke}\rm (i) We refer to \cite[Corollary 2.8]{bauschke-wang} for conditions which ensure the inequality \eqref{pleq}. Let us notice that in case $n=2$, the inequality \eqref{pleq} is automatically fulfilled (cf. \cite[Proposition 3.2(v)]{voisei}).

(ii) In the particular case when $X$ is a Hilbert space and $B:X\rightarrow X$ is the identity operator, the above theorem becomes the extension result proved by Bauschke and Wang in \cite[Theorem 3.2]{bauschke-wang}.
\end{remark}

In the above proof the strong coercivity of the operator $B$ delivers the continuity of the function $f^*$ on the whole space $X^*$ and thus the Fenchel duality result is applicable. Let us note that in order to apply this duality result, it is enough to have a point $(x_0^*,x_0)\in\dom(F_{S,n}^*)$ such that $f^*$ is finite and continuous at $B^*x_0-x_0^*$ (see Theorem \ref{dual-f}). This observation allows us to weaken the strong coercivity of $B$. The price we pay for that is the need to impose a further condition in order to achieve a similar extension result concerning $n$-monotone operators.

\begin{theorem}\label{ext-0-coerc} Let $S:X\rightrightarrows X^*$ be $n$-monotone for some $n\in\{2,3,...\}$. Suppose that $G(S)\neq\emptyset$ and $$p\leq F^*_{S,n}.$$ Consider a linear, continuous, monotone and coercive operator $B:X\rightarrow X^*$. Then for every $w^*\in X^*$ which fulfills the relation \begin{equation}\label{cond-0-coerc}(0,w^*)\in \dom(F_{S,n}^{*\top})-G(B^*)\end{equation} there exists $x\in\cl\co D(S)$ such that $\{(x,w^*-Bx)\}\cup G(S)$ is $n$-monotone.
\end{theorem}

\noindent {\bf Proof.} We follow the lines of the proof of Theorem \ref{ext-1-coerc}. Consider the case $w^*=0$, that is there exists $(x_0^*,x_0)\in\dom(F_{S,n}^*)$ such that $x_0^*=B^*x_0$. With the same notations as above, the function $f$ is coercive, hence $f^*$ is finite and continuous at $0=B^*x_0-x_0^*$. Hence Fenchel duality can be applied (notice that $(RC_1)$ in Theorem \ref{dual-f} is also here fulfilled) and the rest of the proof follows as above.

For the case $w^*\in X^*$ is arbitrary, consider again the operator $S':X\rightrightarrows X^*$ defined by $S'(x)=-w^*+S(x)$ for all $x\in X$. One can prove that $\dom(F_{S',n}^{*\top})=\dom(F_{S,n}^{*\top})-(0,w^*)$ (see Remark \ref{s-tr}) and the conclusion follows.\hfill{$\Box$}

\begin{remark}\label{obs-graph}\rm The condition \eqref{cond-0-coerc} can be replaced by the following one (cf. Proposition \ref{prop-dom})  \begin{equation}\label{cnd-0-coerc-graph}(0,w^*)\in \co G(S)-G(B^*).\end{equation}
\end{remark}

In the above results we have imposed conditions on $B$ in order to obtain some continuity properties of the function $f^*$, which are further used for being able to apply the Fenchel duality result. In the following we would like to notice that one can also use the interior-type regularity conditions in order to ensure the strong duality result in \eqref{dualitate}. In this case, instead of the coercivity of $B$, we guarantee that $$0 \in \core(\dom f^*+A^*(\dom F_{S,n}^*))$$ (the notations are the ones from the proof of Theorem \ref{ext-1-coerc}). Let us define the function $h:X\rightarrow \R$, by $h(x)=\langle Bx,x\rangle$ for all $x\in X$, where $B:X\rightarrow X^*$ is a linear continuous and monotone operator. Let us notice that under this hypotheses $$\dom f^*=\dom h^*+\dom\sigma_C = \dom h^* + \ba C$$ (since $f^*=h^*\Box\sigma_C$, see \cite[Theorem 2.8.7]{Zal-carte}) and 
$$\im(B+B^*)\subseteq \dom h^*.$$
We obtain the following extension result (the details of the proof rely on applying Theorem \ref{dual-f}; thus one has to guarantee $(RC_2)$ in order to get \eqref{dualitate}).

\begin{theorem}\label{ext-core} Let $S:X\rightrightarrows X^*$ be $n$-monotone for some $n\in\{2,3,...\}$. Suppose that $G(S)\neq\emptyset$ and $$p\leq F^*_{S,n}.$$
Consider a linear, continuous and monotone operator $B:X\rightarrow X^*$. Then for every
\begin{equation}\label{cond-core}w^*\in \core\big(\{x^*-B^*x:(x^*,x)\in\dom F_{S,n}^*\}+\im(B+B^*)+\ba D(S)\big)\end{equation} there exists $x\in \cl \co D(S)$ such that $\{(x,w^*-Bx)\}\cup G(S)$ is $n$-monotone.
\end{theorem}

\begin{remark}\label{mult-n-geq-e} \rm (i) By making use of Proposition \ref{prop-dom}, in the condition \eqref{cond-core} one can write $\co G(S)$ in place of $\dom F_{S,n}^{*\top}$ and the extension theorem remains valid. In view of the same result, for $n\geq 3$, we have $\co R(S)\times\co D(S)\subseteq \dom F_{S,n}^*$ and one gets a similar statement for those
 \begin{equation}\label{cond-core2}w^*\in \core\big(\co R(S)-B^*(\co D(S))+\im(B+B^*)+\ba D(S)\big).\end{equation}

(ii) Different to Theorem \ref{ext-1-coerc} and Theorem \ref{ext-0-coerc} the above theorem allows the formulation of an extension result even if $B=0$. In this situation, the conditions \eqref{cond-core} and \eqref{cond-core2} become \begin{equation}\label{cond-core-B0}w^*\in\core\big({\pr}_{X^*}(\dom F_{S,n}^*)+\ba D(S)\big)\end{equation}  and, respectively, \begin{equation}\label{cond-core2-B0}w^*\in\core\big(\co R(S)+\ba D(S)\big).\end{equation} Finally, let us observe that \eqref{cond-core2-B0} implies \eqref{cond-core-B0} (cf. Proposition \ref{prop-dom}). This means that in case $B=0$, for all $w^* \in \core\big(\co R(S)+\ba D(S)\big)$ there exists $x\in \cl \co D(S)$ such that $\{(x,w^*)\}\cup G(S)$ is $n$-monotone.
\end{remark}

\section{Conclusions and further research}

We give in this paper a positive answer to Bauschke and Wang's question (see \cite[Remark 3.5(6)]{bauschke-wang}) concerning whether the convex-analytical approach they propose to obtain extension results for $n$-monotone operators can be extended to the framework of reflexive Banach spaces. We obtain three extension results that rely on the same technique.

We remark that the linear, continuous and monotone operator $B:X\rightarrow X^*$ was used in order to make this extension possible. It could be a topic for further research trying to find out
if the same technique can be adapted to the case $B:X\rightrightarrows X^*$ is a monotone linear relation. Let us recall that $B$ is said to be linear relation if $G(B)$ is a linear subspace of $X\times X^*$. The adjoint of $B$, also denoted by $B^*$, is defined by $$G(B^*)=\{(x,x^*)\in X\times X^*:(x^*,-x)\in (G(B))^{\perp}\},$$ where for any subset $C$ of a topological vector space $Y$ with continuous dual space $Y^*$, $C^{\perp}$ is the \emph{annihilator} of $C$, defined as usual by $C^{\perp}:=\{y^*\in Y^*:\langle y^*,c\rangle=0 \ \forall c\in C\}$. In case $B$ is a monotone linear relation the function  $q_B:X\rightarrow\B$,
$$q_B(x)=\left\{
\begin{array}{ll}
\frac{1}{2}\langle Bx, x\rangle, & \mbox {if } x\in D(B),\\
+\infty, & \mbox{otherwise},
\end{array}\right.$$
is single-valued and convex (cf. \cite[Proposition 2.3]{byw}). It would be interesting to know if $q_B$ can be used instead of $h$ (see the proof of Theorem \ref{ext-1-coerc}) in order to obtain similar extension results for $n$-monotone operators, this time when $B$ is a monotone linear relation. For more on monotone linear relations we refer to \cite{byw, bbw}.

\end{document}